\documentclass[12pt]{amsart}
\usepackage{amssymb}
\newcommand{\query}[1]%
{\mbox{}\marginpar{\raggedright\hspace{0pt}{\small\em #1}}}%
\theoremstyle{plain}

\newtheorem{thm}{Theorem}
\newtheorem{cor}{Corollary}
\newtheorem{lem}{Lemma}
\newtheorem{prop}{Proposition}

\theoremstyle{definition}
\newtheorem{defi}{Definition}
\newtheorem{conj}{Conjecture}
\newtheorem{conv}{Convention}
\newtheorem{nota}[conv]{Notation}
\newtheorem{rem}[thm]{Remark}
\newtheorem{rems}[thm]{Remarks}
\newtheorem{exa}{Example}
\newtheorem{exas}{Examples}
\newtheorem*{sit}{}

\newcommand{\brem}{\begin{rem}}
\newcommand{\brems}{\begin{rems}}
\newcommand{\erem}{\end{rem}}
\newcommand{\erems}{\end{rems}}
\newcommand{\bexa}{\begin{exa}}
\newcommand{\bexas}{\begin{exas}}
\newcommand{\eexa}{\end{exa}}
\newcommand{\eexas}{\end{exas}}
\newcommand{\bdefi}{\begin{defi}}
\newcommand{\edefi}{\end{defi}}
\newcommand{\bcor}{\begin{cor}}
\newcommand{\ecor}{\end{cor}}
\newcommand{\blem}{\begin{lem}}
\newcommand{\elem}{\end{lem}}
\newcommand{\bconv}{\begin{conv}}
\newcommand{\econv}{\end{conv}}
\newcommand{\bconj}{\begin{conj}}
\newcommand{\econj}{\end{conj}}
\newcommand{\bprop}{\begin{prop}}
\newcommand{\eprop}{\end{prop}}
\newcommand{\bthm}{\begin{thm}}
\newcommand{\ethm}{\end{thm}}
\newcommand{\bnota}{\begin{nota}}
\newcommand{\enota}{\end{nota}}
\newcommand{\bsit}{\begin{sit}}
\newcommand{\esit}{\end{sit}}
\newcommand{\be}{\begin{eqnarray}}
\newcommand{\ee}{\end{eqnarray}}

\newcommand{\Spec}{\operatorname{Spec}}

\newcommand{\Aut}{{\operatorname{Aut}}}
\newcommand{\Jac}{{\operatorname{jac\,}}}
\newcommand{\SL}{{\operatorname{SL}}}
\newcommand{\GL}{{\operatorname{GL}}}

\newcommand{\A}{{\mathbb A}}

\newcommand{\C}{{\mathbb C}}

\newcommand{\Z}{{\mathbb Z}}
\newcommand{\N}{{\mathbb N}}

\title{V\'en\'ereau polynomials and related fiber bundles}

\author{Shulim Kaliman}
\address{Department of Mathematics,
University of Miami,
Coral Gables, FL  33124, U.S.A.}
\email{Shulim.Kaliman@math.miami.edu}

\author{Mikhail Zaidenberg}
\address{Universit\'e
Grenoble I, Institut Fourier, UMR 5582 CNRS-UJF, BP 74,
38402 St.\ Martin
d'H\`eres cedex, France}
\email{Mikhail.Zaidenberg@ujf-grenoble.fr}

\thanks{
{\bf Acknowledgments:} This research
started during a visit of the first
author at the Institut Fourier of
the University of Grenoble, and continued
during a stay of the second author
at the Max Planck Institute of Mathematics
at Bonn.
The authors thank both institutions for
their support.
It is our pleasure to thank
Don Zagier for stimulating discussions,
and M. Uludag for his help with MAPLE}

\thanks{
\mbox{\hspace{11pt}}{\it 1991 Mathematics
Subject Classification}:
14R10, 14R25.\\
\mbox{\hspace{11pt}}{\it Key words}:
polynomial ring, variable,
algebraic fiber bundle, flat family.}

\begin{document}

\begin{abstract}
The V\'en\'ereau polynomials
$$v_n:=y+x^n(xz+y(yu+z^2)),\qquad n\ge 1,$$
on $\A_\C^4$ have all fibers isomorphic to the affine space
$\A_\C^3$. Moreover, for all $n\ge 1$ the map $(v_n, x):\A_\C^4\to
\A_\C^2$ yields a flat family of affine planes over $\A_\C^2$. In
the present note we show that over the punctured plane
$\A_\C^2\backslash \{\bar 0\}$, this family is a fiber bundle.
This bundle is trivial if and only if $v_n$ is a variable of the
ring $\C[x][y,z,u]$ over $\C[x]$.

It is an open question whether $v_1$ and $v_2$ are variables of
the polynomial ring $\C^{[4]}=\C[x,y,z,u]$, whereas S.
V\'en\'ereau established that $\,v_n$ is indeed a variable of
$\C[x][y,z,u]$ over $\C[x]$ for $n\ge 3$. In this note we give
another proof of V\'en\'ereau's result based on the above
equivalence. We also discuss some other equivalent properties, as
well as the relations to the Abhyankar-Sathaye Embedding Problem
and to the Dolgachev-Weisfeiler Conjecture on triviality of flat
families with fibers affine spaces.
\end{abstract}

\maketitle

In \cite{KVZ1, KVZ2} polynomials in four variables of the form
\be\label{gen} p=f(x,y)u+g(x,y,z)\,\ee were studied. It was shown
that $p:\A_\C^4\to \A_\C^1$ is a flat family of affine spaces
(i.e., every fiber $p^*(c)$, $c\in \A_\C^1$, is reduced and
isomorphic to $\A_\C^3$) provided that $p^*(0)\simeq \A_\C^3$
(\cite[Theorem 3.21]{KVZ2}). As for the latter condition, a
criterion is given in terms of the morphism $\pi:C\to\Gamma$,
where $C=\{f=g=0\}\subseteq \A_\C^3$ and $\Gamma=\{f=0\}\subseteq
\A_\C^2$ are affine curves and $\pi$ is the projection
$(x,y,z)\longmapsto (x,y)$ \cite[Theorems 2.11 and 3.21]{KVZ2}.
Based on this, S. V\'en\'ereau \cite{Ve} considered the following
polynomials generating flat families of affine spaces:
\be\label{ven} v_n:=y+x^n(xz+y(yu+z^2))
=x^ny^2u+y+x^{n+1}z+x^nyz^2,\qquad n\ge 1\,.\ee He showed that for
every $n\ge 3$, $v_n$ is a variable of the ring $\C^{[4]}$, or in
other words, a coordinate of a polynomial automorphism of
$\A_\C^4$. Furthermore, some sufficient conditions were found (see
\cite[Section 4.1]{KVZ2}) for a polynomial $p$ of form (\ref{gen})
to be a variable. However we do not know whether the latter is
true if these sufficient conditions are not satisfied, and this is
the case for $v_1$ and $v_2$.  That is, we do not know whether the
hypersurfaces $v_i^{-1}(0)\simeq \A_\C^3$, $i=1,2$, in $\A_\C^4$
can be rectified by means of polynomial automorphisms of $\A_\C^4$
(this is a particular case of the Abhyankar-Sathaye Embedding
Problem).

It was also shown in \cite{Ve} that
for every $n\ge 3$ the morphism
$$\Phi_n=(v_n, x):\A_\C^4\to \A_\C^2$$ defines
a trivial family of affine planes.
This  means that the V\'en\'ereau polynomials
$v_n$ ($n\ge 3$) are
$x$-variables i.e., variables
of the $\C[x]$-algebra $\C[x][y,z,u]$.

If a polynomial $p$ as in (\ref{gen})
defines a flat family of affine
spaces $p:\A_\C^4\to \A_\C^1$ then $p$ also
defines, along with a suitable second polynomial
$q\in\C[x,y]\subseteq \C^{[4]}$, which is
a variable in $\C[x,y]$,
a flat family of affine planes
$\Phi=(p,q): \A_\C^4\to \A_\C^2$
(see \cite[Theorem 3.21]{KVZ2}).
Thus the question arises whether
the two remaining families of affine planes
$\Phi_n=(v_n, x):\A_\C^4\to \A_\C^2$ ($n=1,2$)
and, more generally,
all such families of affine planes
$\Phi=(p,q): \A_\C^4\to \A_\C^2$
as above,
are trivial.

In the present note we address the former special question
by showing in Proposition \ref{equi} that
$\Phi_n$ restricts
to an algebraic fiber bundle, say, 
$\lambda_n$ over the punctured plane:
$$\lambda_n:=\Phi_n|(\A_\C^4\backslash \A_\C^2):
\A_\C^4\backslash \A_\C^2\to
\A_\C^2\backslash \{\bar 0\}\,,$$
where $\A_\C^2\subseteq \A_\C^4$ is the coordinate
$(z,u)$-plane.
This bundle has $\A_\C^2$ as
the typical fiber and
$\Aut (\A_\C^2)$ as the structure group.
 We show that
$\lambda_n$ is trivial for every $n\ge 3$, thus recovering
V\'en\'ereau's result (see Corollary \ref{vcor} below).
In the cases $n=1,2$
we were not able to carry the
computations needed by our methods,
and so the question remains open.

Generally speaking, we deal below with the following
three categories of affine
$S$-schemes $f:X\to S$ over a quasiprojective base $S$:

\smallskip

$\bullet$ The flat families of affine $m$-spaces over $S$
\footnote{That is, every fiber $X_s=f^*(s),\,\,s\in S$,
is a reduced scheme isomorphic to $\A_\C^m$.};

\smallskip

$\bullet$ The algebraic fiber bundles over $S$
with fiber $\A_\C^m$;

\smallskip

$\bullet$ The algebraic vector bundles of rank
$m$ over $S$.

\smallskip

\noindent
The Dolgachev-Weisfeiler Conjecture \cite[(3.8.3)]{DW}
claims that the first category reduces to the second one.
In turn, for an affine base $S$,
the second one  reduces to the third one
due to the
Bass-Connell-Wright Theorem \cite{BCW}.

If
the Dolgachev-Weisfeiler Conjecture
were answered in
affirmative this would provide an affirmative
answer to our general question.
Indeed,
in this case $\Phi$ (in particular, $\Phi_n$) would be a fiber
bundle over $\A_\C^2$ with fiber $\A_\C^2$, thus a vector bundle
by the Bass-Connell-Wright Theorem, hence it must be trivial due
to the Quillen-Suslin Theorem.

Conversely, if any one of the bundles
$\lambda_1$ or $\lambda_2$ were non-trivial
this would provide a counterexample
to the Dolgachev-Weisfeiler Conjecture
(with $m=2$), and,
presumably, to
the Abhyankar-Sathaye Embedding Problem.

Summarizing, the following equivalences
hold.

\bprop\label{equi} If $v\in\C^{[4]}$
is a polynomial such that $\Phi=(x,v)$
yields a flat family of affine planes
$\Phi: \A_\C^4\to \A_\C^2$
then the following conditions are equivalent:
\begin{enumerate}
\item[($i$)] $v$ is an $x$-variable of the ring
$\C[x,y,z,u]$.
\item[($ii$)] $\Phi: \A_\C^4\to \A_\C^2$ is an
algebraic fiber bundle.
\item[($iii$)] It is a trivial fiber bundle.
\item[($iv$)] $\lambda:= \Phi|(\A_\C^4\backslash F_{\bar 0})$,
where $F_{\bar 0}= \Phi^{-1}(\bar 0)$, is a trivial
fiber bundle over
$\A_\C^2\backslash \{\bar 0\}$.
\end{enumerate}\eprop

\proof The equivalence ($i$)$\Leftrightarrow$($iii$)
is a tautology, and
($iii$)$\Leftrightarrow$($ii$)
follows by combining the Bass-Connell-Wright
Reduction Theorem
and the Quillen-Suslin Theorem as above.
The implication ($iii$)$\Rightarrow$($iv$) is evident,
whereas the converse one
($iv$)$\Rightarrow$($iii$) can be easily obtained
by extending
the trivialization morphism
$$\A_\C^4\backslash F_{\bar 0}\to
(\A_\C^2\backslash \{\bar 0\})\times \A_\C^2$$ and its inverse to
the deleted planes $F_{\bar 0}$ and $\{\bar 0\}\times \A_\C^2$,
respectively. \qed\medskip

\bnota\label{nota} We use below the following polynomials
from the ring $\C[x,y,z,u]$ :
\be\label{w}   w & =   z^2+yu\,,\\
\label{t}   t   & =   xz+yw\,,\\
\label{s}  s  & =-(2xz+yw)w\,,\\
\label{eta}   \eta  &  =s+x^2u=x^2u-2xzw-yw^2\,,\\
\label{zeta1}   \zeta^{(1)} &
 =-v_1z+v_1t(v_1u+z^2+w)+t^2(xz^2+st)\,,\\
\label{zeta2}   \zeta^{(2)} &
  =-z+xt(v_2u+z^2+st+w)\,,\\
\label{zetan}   \zeta^{(n)} &
 =-z+x^{n-3}t\left(v_n\eta+x^2w\right),
\qquad n\ge 3\,,\ee
where
\be\label{vn} v_n=y+x^nt\ee
stands for
the V\'en\'ereau polynomials (\ref{ven}).
We also consider the rational functions
\be\label{rafa} \xi_0=\eta/x^3,\qquad\xi_1^{(1)}
=\zeta^{(1)}/v_1^3\qquad
\mbox{and}\qquad \xi_1^{(n)}
=\zeta^{(n)}/v_n^2, \quad n\ge 2\,.\ee
We let
$$L_0=\C[x,x^{-1}],\quad
L_1=\C[v,v^{-1}],\quad
K_0=\\C[x,x^{-1},v],\quad K_1=\C[x,v,v^{-1}]\,,$$
$$R=\C[x,v]=K_0\cap K_1\qquad
\mbox{and}\qquad M= \C[x,x^{-1},v,v^{-1}]\,,$$  with the
convention that $v=v_n$ whenever we consider $\Phi_n$ or
$\lambda_n$. Letting $\A_\C^2=\Spec\C[x,v_n]$, the punctured plane
$\A_\C^2\backslash \{\bar 0\}$ can be covered by the Zariski open
subsets
$$U_0=\Spec K_0=\A_\C^2\backslash \{x=0\}
\qquad\mbox{and}\qquad
U_1=\Spec K_1=\A_\C^2\backslash \{v_n=0\}\,,$$
where $U_0\cap U_1=\Spec M$. \enota

With this notation the following results hold.

\bprop\label{loctri} (a) The morphism
$\varphi_0=(t,\xi_0)$,
respectively,
$\varphi_1^{(n)}=\left(t,\xi_1^{(n)}\right)$,
yields a local
trivialization for the family $\lambda_n$
over $U_0$, respectively, over $U_1$.
Thus $\lambda_n$ is an algebraic fiber bundle
with $\A_\C^2$ as the typical fiber.

(b)
The transition function
$$\varphi_{10}^{(n)}=\varphi_1^{(n)}\circ
\varphi_0^{-1}:(t,\xi_0)\longmapsto \left(t, \xi_1^{(n)}\right)$$
is a triangular automorphism with
\be\label{trfun} \xi_1^{(n)}=\xi_0
+\frac{p_n(t)}{x^kv_n^l},
\qquad\mbox{where}\qquad
(k,l)=\begin{cases}
(3,2),\qquad & n\ge 2\\
(3,3),\qquad & n=1\,,
\end{cases}
\ee
and the polynomials $p_n\in\C[x,v_n,t]$
are given by
\be\label{pn}
p_n(t)=
\begin{cases}
{v_nt^2-x^2t,}\qquad & {n\ge 3}\\
{x^2t^3+v_2t^2-x^2t,}\qquad & {n=2}\\
{x^2t^4+xv_1t^3+v_1^2t^2-x^2v_1t,}\qquad & {n=1}\,.
\end{cases}
\ee
\eprop

\proof (a) The Nagata automorphism $\alpha\in \Aut (\A_\C^3)$ of
$\A_\C^3=\Spec \C[y,z,u]$ is given (see \cite{Na}) by
 $\alpha: (y,z,u)\to (y,t_0,\eta_0)$, where
$$t_0=z+yw\qquad \mbox{and}\qquad \eta_0=u-2zw-yw^2\,$$
with $w=z^2+yu$ as in (\ref{w}) above. 
\footnote{Let us observe in
passing that $\alpha=\mu\circ\delta\circ\mu^{-1}$, 
where $\mu:
(y,z_1,u_1)\longmapsto (y,z=yz_1,u=yu_1)$, or else $\mu:
(y,t_1,\eta_1)\longmapsto (y,t_0=yt_1,\eta_0=y\eta_1)$ is a
birational morphism, or in other words, 
an affine modification of
$\A_\C^3$ with the locus $(D,\bar 0)$, 
where $D=\{y=0\}$, and
$\delta:=\delta_3\circ\delta_2\circ\delta_1\in 
\Aut (\A_\C^3)$ is a
tame $\C[y]$-automorphism with 
$$\delta_1:(y,z_1,u_1)\longmapsto
(y,z_1,w_1)\,,\qquad w_1:=u_1+z_1^2,$$
$$\delta_2: (y,z_1,w_1)\longmapsto
(y,t_1,w_1)\,,\qquad t_1:=z_1+y^2w_1,$$
$$\delta_3: (y,t_1,w_1)\longmapsto
(y,t_1,\eta_1)\,,\qquad \eta_1:=w_1-t^2_1\,.$$}
Composing $\alpha$ with the following
$L_0$-automorphisms
of $L_0[y,z,u]$ (cf. \cite{Ve}):
$$g: (y,z,u)\longmapsto (x^{-2}y,xz,x^4u)
\qquad\mbox{and}\qquad
h:(y,z,u)\longmapsto (x^2y,z,x^{-2}u)\,$$
we obtain:
$$w\circ g=x^2w,\qquad t_0\circ g=t,\qquad
\eta_0\circ g=x^2\eta\,$$
and $$\beta:=h\circ\alpha\circ g
\in \Aut_{L_0} (L_0[y,z,u]),\qquad\mbox{where}
\qquad  \beta:
(y,z,u)\longmapsto (y,t,\eta)\,.$$
Letting
$$\gamma_n:(y,z,u)\longmapsto (v_n=y+x^nt,t,\eta)$$
it follows that $\gamma_n\circ\beta\in\Aut_{L_0}
(L_0[y,z,u])$.
Thus
$$\Phi_n^{-1} (U_0)=\A_\C^4\backslash \{x=0\}=
\Spec\C[x,x^{-1},y,z,u]=\Spec\, \C[x,x^{-1},v_n,t,\eta] \cong
U_0\times\A_\C^2$$ with $\A_\C^4=\Spec\,\C[x,y,z,u]$ and
$\A_\C^2=\Spec\C[t,\eta]$. Hence
$\varphi_0=(t,\xi_0)=(t,\eta/x^3)$ yields indeed a trivialization
of $\lambda_n$ over $U_0$.

To show that also $\varphi_1^{(n)}=(t,\xi_1^{(n)})$
yields
a trivialization of $\lambda_n$ over $U_1$,
we consider
separately the cases $n=1,\,n=2$ and $n\ge 3$.
We will constantly exploit the relations
\be\label{rel1}  ys+t^2=x^2z^2
\qquad\Rightarrow\qquad v_ns+t^2=x^2z^2+x^nst\,\ee
\be\label{rel2}
\quad\Rightarrow\quad v_n\eta+t^2=x^2(z^2+v_nu)+x^nst\ee
(see (\ref{eta}) and (\ref{vn})).

\smallskip

{\it Case $n=1$.} We denote
\be\label{zeta'} \zeta':=\frac{v_1\eta+t^2}{x}
=x(v_1u+z^2)+st\ee
(see (\ref{rel2})),
\be\label{zeta''} \zeta'':  =\frac{v_1\zeta'+t^3}{x}=
v_1(v_1u+z^2)+xz^2t+st^2\ee
$$ =yw+xtw+xtv_1u+xz^2t+st^2$$
(see (\ref{w}) and (\ref{vn})), and
\be\label{zeta'''}\zeta''':=\zeta''-t
=-xz+xtw+xtv_1u+xz^2t+st^2\ee
(see (\ref{t})).
Now
\be\label{zetarel} \zeta^{(1)}:
=\frac{v_1\zeta'''+t^4}{x}=
\frac{v_1^3\eta+p_1(t)}{x^3}\,\ee
verifies both (\ref{zeta1})
and (\ref{trfun}) for $n=1$.

Further, for any point $C=(c_1,c_2)\in \A_\C^2$ with $c_1\neq
0,\,c_2\neq 0$, the functions $(t,\eta)$, hence also
$(t,\zeta^{(1)})$, give global coordinates on the fiber
$F_{C}=\Phi_1^{-1}(C)\simeq\A_\C^2$ over $C$. In the case
$c_1=0,\,c_2\neq 0$ the following hold: \be\label{rela0}
v_1=y=c_2,\quad w=c_2^{-1}t=c_2u+z^2, \quad
\eta=s=-c_2^{-1}t^2\,,\ee and \be\label{rela1}  \zeta^{(1)}=
-c_2z+2t^2-c_2^{-1}t^5\,\ee
 (see (\ref{w}), (\ref{zeta1}), (\ref{vn})
and (\ref{rel1})).
Therefore we obtain:
\be\label{iso}  \C[z,u]=\C[z,t]=
\C[t,\zeta^{(1)}]=\C[t,\xi^{(1)}_1]\,. \ee
 Clearly, $F_C\simeq\A_\C^2=\Spec\C[z,u]$
for $C=(0,c_2)$
with $c_2\neq 0$. In other words,  $(z,u)$ give coordinates
on the fiber $F_C$. Thus by (\ref{iso}) for $c_2\neq 0$
the
functions $(t,\zeta^{(1)})$,
and hence also $(t,\xi_1^{(1)})$,
provide as well coordinates on $F_C$.
Now (a) follows for $n=1$.

\smallskip

{\it Case $n=2$.}  From (\ref{rel2})
we obtain:
\be\label{zeta2'} \zeta':=
\frac{v_2\eta+t^2}{x^2}=z^2+uv_2+st\,.\ee
It can be easily seen that the polynomial
\be\label{zeta2''} \zeta^{(2)}
:= \frac{v_2\zeta'+t^3-t}{x}=
\frac{v_2^2\eta+p_2(t)}{x^3}\ee
verifies both (\ref{zeta2}) and (\ref{trfun})
for $n=2$. It follows that
for $C=(c_1,c_2)$ with $c_1\neq 0, \,c_2\neq 0$,
the map $(t,\eta)\longmapsto (t,\zeta^{(2)})$
provides an isomorphism
$$F_C\simeq\A_\C^2=\Spec\,\C[t,\eta]=
 \Spec\C[t,\zeta^{(2)}]\,.$$
As above,
for $c_1=0,\,c_2\neq 0$ we get
$$ \C[z,u]=\C[z,t]=\C[t,\zeta^{(2)}]\,$$
and so, the functions
$(t,\zeta^{(2)})$
still give coordinates on the fiber $F_C$.
This proves (a) for $n=2$.

\smallskip

{\it Case $n\ge 3$.} Letting again
 \be\label{zetan'} \zeta':=\frac{v_n\eta+t^2}{x^2}
=z^2+v_nu+x^{n-2}st\,\ee
(see (\ref{rel2})) we can easily see that
\be\label{zetan''} \zeta^{(n)}:
=\frac{v_n\zeta'-t}{x}=
\frac{v_n^2\eta+p_n(t)}{x^3}\ee
verifies both
(\ref{zetan}) and (\ref{trfun}), the map
$(t,\eta)\longmapsto (t,\zeta^{(n)})$ provides
an isomorphism $ F_C\simeq \A_\C^2=\Spec
\C[t,\zeta^{(n)}]$
as soon as $c_1\neq 0$
and $c_2\neq 0$,
and, moreover, the functions $(t,\zeta^{(n)})$
yield coordinates
on any fiber $F_C$ with $c_2\neq 0$.
This proves (a) in  the  general case.
Actually (b) has been established
in the course of  proof of (a).
\qed\medskip

\brem\label{dejo} By Proposition \ref{loctri}(b)
the transition function $\varphi_{10}^{(n)}$
of $\lambda^{(n)}$ as in (\ref{trfun})
takes values in the subgroup of plane
triangular automorphisms of the form
$$(t,\xi)\longmapsto
\left(t,\xi+\frac{q(t)}{x^nv^m}\right)
\qquad\mbox{with}\qquad q\in R[t],\quad
n,m\in\N\cup\{0\}\,.$$
However, within this subgroup, the fiber bundles
$\lambda_n$ are non-trivial.
Indeed, triviality of $\lambda_n$ is equivalent
to the existence of
a decomposition $$ \frac{p_n(t)}{x^kv^l}
=\frac{a_n(t)}{x^{k+\alpha}}
-\frac{b_n(t)}{v^{l+\beta}}\qquad\mbox{with}\qquad
a_n,\,b_n\in R[t],\quad\alpha,\beta\ge 0\,,$$
and hence is equivalent to:
$p_n\in(x^k,v^l)R[t]$,
which is not the case (see (\ref{pn})).
\erem\medskip

The following results will be useful for establishing the
triviality of bundles $\lambda_n$.

\bprop\label{triv} Let $\lambda$ be an
algebraic fiber bundle
over $\A_\C^2\backslash \{\bar 0\}$,
where $\A_\C^2=\Spec\C[x,v]$,
with $\A_\C^2=\Spec\C[t,\xi]$ as the typical fiber.
We suppose that the restrictions $\lambda|U_i$
($i=0,1$) are trivial. If the transition function
$$\varphi_{10} : U_0\cap U_1\to\Aut(\C[t,\xi])$$
has the form:
\be\label{form} \varphi_{10}:(t,\xi)\longmapsto
\left(t,\xi+\frac{p(t)}{x^kv^l}\right)\ee
with $p\in R[t]$, where $R=\C[x,v]$ and $k,l\in\N$,
then the following hold.

(a)
$\lambda$ is trivial if and only if
\be\label{tri} \varphi_{10}
=\tau_1\circ\tau_0^{-1}
\qquad\mbox{with}\qquad
\tau_i\in \Aut_{K_i} K_i[t,\xi],\,\,i=0,1\,.\ee
If (\ref{tri}) is fulfilled then necessarily
$\tau_0$ and
$\tau_1$ have the form:
\be\label{tau} \tau_0:(t,\xi)\longmapsto
\left(a,\,\,\frac{b_0}{x^k}\right)
\qquad\mbox{and}\qquad
\tau_1:(t,\xi)\longmapsto
\left(a,\,\,\frac{b_1}{v^l}\right)\,,\ee
where $a,\,b_0,\,b_1\in R[t,\xi]$ satisfy the cocycle relation
\be\label{relp} x^kb_1-v^lb_0=p(a)\,.\ee
Furthermore, up to multiplying $b_0,\,b_1$
by a constant
we may assume that the following Jacobian relations hold:
\be\label{jac} \Jac (\tau_0)=1\Leftrightarrow
\Jac (a,b_0)=x^k,
\qquad
\Jac (\tau_1)=1\Leftrightarrow
\Jac (a,b_1)=v^l\,,\ee
where $\Jac (*,*)$ stands for the
jacobian in $t,\xi$.

\smallskip

(b)
Any solution $(a,\,b_0,\,b_1)$ in
$(R[t,\xi])^3$ of (\ref{relp})
satisfies
\be\label{d} \frac{\Jac (a,b_0)}{x^k} =
\frac{\Jac (a,b_1)}{v^l}=:
d\in R\ee
and
\be\label{prim} \Jac(b_0,b_1)
=-d\cdot p'(a)\,.\ee
Such a solution also verifies (\ref{jac})
if and only if $d=1$.
In the latter case
\be\label{newjac} \Jac (b_0,b_1)=-p'_t(a)\,.\ee
\eprop

\proof Although (a) is well known, we remind the proof. We have:
$\varphi_{10}= \varphi_1\circ\varphi_0^{-1}$, where $\varphi_i$ is
a trivialization of $\lambda|U_i$ ($i=0,1$). If (\ref{tri}) holds
then $\psi|U_i:=\tau_i^{-1} \varphi_i$ ($i=0,1$) defines a global
trivialization of $\lambda$. Indeed, $\psi$ is a well defined
regular map over $U_0\cup U_1=\A_\C^2\backslash\{\bar 0\}$, since
by (\ref{tri}), $\tau_0^{-1}\varphi_0=\tau_1^{-1}\varphi_1$ in
$U_0\cap U_1$. Conversely, if $\psi$ is a global trivialization of
$\lambda$ then $\tau_i:=\varphi_i\circ\psi^{-1}\in \Aut_{K_i}
K_i[t,\xi]$ ($i=0,1$) satisfy (\ref{tri}), as required.

Letting $\tau_i=(\hat a_i,\,\hat b_i)$ with $\hat a_i,\,\hat
b_i\in K_i[t,\xi]$, by (\ref{tri}) we obtain: \be\label{triv1}
(\hat a_1,\,\hat b_1) =\varphi_{10}\circ (\hat a_0,\,\hat b_0)
=\left(\hat a_0,\,\hat b_0+ \frac{p(\hat
a_0)}{x^kv^l}\right)\,.\ee Therefore $\hat a_0=\hat a_1=:a\in
R[t,\xi]$, and \be\label{b} \hat b_0
=\frac{b_0}{x^{k+\alpha}},\qquad \hat b_1
=\frac{b_1}{v^{l+\beta}}\,,\ee where $b_0,\,b_1\in R[t,\xi]$ and
$\alpha,\,\beta\ge 0$. Now (\ref{b}) yields: \be\label{rl}
x^{k+\alpha} b_1-v^{l+\beta}b_0 =x^{\alpha}v^{\beta}p(a)\,,\ee
hence $x^{\alpha}|b_0$ and $v^{\beta}|b_1$. Thus without loss of
generality we can suppose that $\alpha=\beta=0$, and so
(\ref{relp}) follows.

As $\Jac (\tau_i)$ is a
unit in $K_i[t,\xi]$  ($i=0,1$) we have
$$\Jac (\tau_0)=c_0x^n\quad\mbox{and}\quad
\Jac (\tau_1)=c_0v^m
\qquad\mbox{for some}\qquad
n,m\in\Z \quad\mbox{and}\quad
c_0,c_1\in\C\backslash \{0\}\,.$$
Since
$\Jac (\varphi_{10}) = 1$
it follows from (\ref{tri}) that
$$\Jac (\tau_0)
=\Jac (\tau_1)\in K_0\cap K_1=R\,$$ is a unit in $R$. Hence
$n=m=0$, and we may suppose that $c_0=c_1=1$, which yields
(\ref{jac}).

\noindent (b) Applying to (\ref{relp})
the derivations $*\longmapsto \Jac (a,*)$
and
$*\longmapsto \Jac (*, b_1)$, respectively,
gives
\be\label{rel0} x^k\Jac (a,b_1)=v^l
\Jac (a,b_0)\qquad\mbox{and}\qquad
p'(a)\Jac (a,b_1)=-v^l\Jac (b_0,b_1)\,.\ee
Now (\ref{d}) and (\ref{prim}) follow.
\qed\medskip

\brems \label{remk} 1. If $(a,b_0,b_1)$
is a solution of (\ref{relp}) then so is
$(a,\,b_0+x^kc,\,b_1+v^lc)$ for any $c\in R[t,\xi]$.
This new solution verifies (\ref{jac}) if and only if
$\Jac (a,c)=1-d$ with $d$ as in (\ref{d}).

Similarly, $(a+A,\,b_0+B_0,\, b_1+B_1)$
with $A,B_0,B_1\in R[t,\xi]$ gives
a new solution of
(\ref{relp})
if and only if
\be\label{tailor} \sum_{j=1}^{\deg p}
\frac{p^{(j)}(a)}{j!}A^j=x^kB_1-v^lB_0\,.\ee

2. By virtue of (\ref{d}), the cocycle relation (\ref{relp}) and one of
the Jacobian relations (\ref{jac}) imply the other one.

3. Modulo the plane Jacobian Conjecture the maps $\tau_0$ and
$\tau_1$ as in (\ref{tau}) are invertible if and only if
the Jacobian relations (\ref{jac}) hold. 
Presumably, in our particular setting the
degrees of the coordinate polynomials 
can be limited to the range where the
plane Jacobian Conjecture is known to be true. 
If so then the triviality of the bundles $\lambda_n$ 
reduces to the existence of a solution 
$(a,b_0,b_1)\in (R[t,\xi])^3$ of 
both (\ref{relp}) and (\ref{jac}).

\erems

\blem\label{l1} Suppose that $p=p_n$ in (\ref{relp}) and
$(k,l)=(3,2)$ if $n\ge 2$, $(k,l)=(3,3)$ if $n=1$. Then there
exists a solution $(a,\,b_0,\,b_1)\in (R[t,\xi])^3$ of
(\ref{relp}) with \be\label{ij} a=\sum_{i,j\ge 0} a_{ij} x^iv^j,
\qquad\mbox{where}\qquad a_{ij}\in S: =\C[t,\xi]\,,\ee if and only
if \be\label{cond} a_{00}=0\qquad\mbox{and}\qquad
a_{01}=a_{10}^2\,.\ee \elem

\proof Clearly, (\ref{relp}) has a solution
if and only if
\be\label{con}
p(a)=p_n(a)\in
(x^k,v^l),\qquad\mbox{where}\qquad (k,l)=
\begin{cases}
(3,2), & n\ge 2\\
(3,3), & n=1\,.
\end{cases}
\ee
 We have
\be\label{co}  p_n(a)\equiv
\begin{cases}
a_{00}^2v \mod (x,v^2),  & n\ge 2\\
a_{00}^3 xv \mod (x^2,v^2),  & n=1
\,.\end{cases}
\ee
In any case
\be\label{c} p_n(a)\in (x^k,v^l)
\Rightarrow a_{00}=0\,,\ee
and then
\be\label{co1}  p_n(a)\equiv
\begin{cases}
(a_{10}^2-a_{01})x^2v \mod (x^3,v^2),  & n\ge 2\\
(a_{10}^2-a_{01})x^2v^2 \mod (x^3,v^3),  & n=1
\,.\end{cases}
\ee
Now the assertion follows. \qed\medskip

\bnota\label{sa} Letting in Proposition
\ref{triv}
$$p(t)=
\sum_{j=0}^
{{\rm deg}\, p} r_jt^j,\qquad\mbox{where}\qquad r_j\in R=\C[x,v]\,,$$
and supposing that $r_0=0$,
we can define
the successive approximations
\be\label{app}
\varphi_{10}^{(m)}:(t,\xi)\longmapsto\left(t,\xi+
\frac{\sum_{j=0}^{m} r_jt^j}{x^kv^l}\right),\qquad
m=1,\ldots,\deg p,\ee
to our transition function $\varphi_{10}$.
These $\varphi_{10}^{(m)}$ may be considered
as the transition functions
of a succession of algebraic fiber bundles
$\lambda^{(m)}$
over the punctured plane
 $\A_\C^2\backslash \{\bar 0\}$
with $\A_\C^2$ as the typical fiber.
In particular, the linear part
$\varphi_{10}^{(1)}$ defines
a rank 2 vector bundle $\lambda^{(1)}$
over the punctured plane. \enota

With these notation the following holds.

\blem\label{l2} If $ r_1=-x^{k'}v^{l'}$, where
$\alpha:=k-k'\ge 0$
and $\beta:=l-l'\ge 0$, then the vector bundle
$\lambda^{(1)}$ is trivial.
In particular, $\lambda^{(1)}_n$ is trivial
for every $n\ge 1$.
\elem

\proof Similarly as in Proposition \ref{triv}(a),
$\lambda^{(1)}$ is trivial if and only if
\be\label{tr} \varphi_{10}^{(1)}
=\tau_1^{(1)}\circ\left(\tau_0^{(1)}\right)^{-1}
\qquad\mbox{with}\qquad
\tau_i^{(1)}\in \GL(2, K_i),\,\,i=0,1\,.\ee
By our assumption the linear part $\varphi_{10}^{(1)}$
of $\varphi_{10}$ is:
$$ \varphi_{10}^{(1)} : (t,\xi)\longmapsto
\left(t,\,\xi- x^{-\alpha}v^{-\beta}t\right)\,.$$
It
can be represented by the following matrix:
\be\label{mat} g=\left(
\begin{matrix}
1 & 0\\
-x^{-\alpha}v^{-\beta} & 1
\end{matrix}\right)\in\SL (2, M)\,.\ee
Letting e.g.,
\be\label{mats} \tau_0^{(1)}=\left(
\begin{matrix}
v^{\beta} & -x^{\alpha}\\
x^{-\alpha} & 0
\end{matrix}\right)\in\SL (2, K_0)
\quad\mbox{and}\quad
\tau_1^{(1)}=\left(
\begin{matrix}
v^{\beta} & -x^{\alpha}\\
0 & v^{-\beta}
\end{matrix}\right)\in\SL (2, K_1)
\,\ee we obtain $g\circ \tau_0^{(1)}=\tau_1^{(1)}$, which shows
that $\lambda^{(1)}$ is indeed trivial (cf. (\ref{tri})).
According to (\ref{pn}) in Proposition \ref{loctri},
$\lambda=\lambda_n$ corresponds to the particular case where
$(\alpha,\beta) =(1,2)$, hence $\lambda^{(1)}_n$ is also trivial.
\qed\medskip

The similar triviality result holds also for the second order
approximations.

\bprop\label{trivmain} The fiber bundle
$\lambda^{(2)}_n$
is trivial for every $n\ge 1$. \eprop

\proof By virtue of (\ref{pn}) all the fiber bundles
$\lambda^{(2)}_n$, $n\ge 1$, actually have the same transition
function \be\label{tf} \varphi_{10}^{(2)}: (t,\xi)\longmapsto
\left(t,\,\,\xi-\frac{t}{xv^2} +\frac{t^2}{x^3v}\right)\,.\ee If a
decomposition \be\label{tr1} \varphi_{10}^{(2)}=
\tau_1^{(2)}\circ\left(\tau_0^{(2)}\right)^{-1}
\qquad\mbox{with}\qquad \tau_i^{(2)}\in \Aut_{K_i}
K_i[t,\xi],\,\,i=0,1\,,\ee as in (\ref{tau}) does exist then
clearly there should also exist such a decomposition with the
linear parts $\tau_i^{(1)}$ of $\tau_i^{(2)}$ ($i=0,1$) as in
(\ref{mats}), where $(\alpha,\beta)=(1,2)$, i.e., with
\be\label{3co} a=v^2t-x\xi+O(2),\qquad
b_0=x^2t+O(2)\quad\mbox{and}\quad b_1=
\begin{cases}
\xi+O(2), & n\ge 2\\
v\xi+O(2), & n=1\,,
\end{cases}
\ee where $O(2)$ stands for the terms 
of order $\ge 2$ in $t,\xi$.
By Lemma \ref{l1}, $a_{00}=0$ and $a_{01}=a_{10}^2$.  This
indicates that the power series expansion of the polynomial
$a\in\C[x,v,t,\xi]$ can start e.g., with: 
\be\label{pi} \tilde a
=v^2t-x\xi+v\xi^2=v\delta-x\xi\,,\ee where 
\be\label{del}
\delta:=vt+\xi^2, \qquad a_{10}=-\xi\quad\mbox{and}\quad
a_{01}=a_{10}^2=\xi^2\,.\ee
 For the fiber bundles $\lambda^{(2)}_n$
the cocycle relation (\ref{relp}) becomes (see (\ref{pn})):
\begin{eqnarray}
-x^2a+va^2 = x^3b_1-v^2b_0,\qquad & n\ge 2 \\
-x^2va+v^2a^2 = x^3b_1-v^3b_0,\qquad &  n=1\,.
\end{eqnarray}
It is satisfied e.g., by the triple $(a,b_0,b_1)$ with
\be\label{sol} a=\tilde a,\qquad
b_0=x^2t-v\delta^2+2x\delta\xi
\qquad\mbox{and}\qquad b_1=
\begin{cases}
\xi, & n\ge 2\\
v\xi, & n=1\,,
\end{cases}\ee
cf. (\ref{3co}) (recall that $(k,l)=(3,2)$ if $n\ge 2$,
$(k,l)=(3,3)$ if $n=1$). Thus for any $n\ge 1$ one can take e.g.,
\be\label{tau02}  \tau_0^{(2)} : (t,\xi)\longmapsto
\left(a,\frac{b_0}{x^k}\right)=
\left(v\delta-x\xi,\,\,\,\frac{1}{x}t -\frac{v}{x^3}\delta^2+
\frac{2}{x^2}\delta\xi\right)\ee and \be\label{tau12} \tau_1^{(2)}
: (t,\xi)\longmapsto\left(a,\frac{b_1}{v^l}\right)=
\left(v\delta-x\xi,\,\,\, \frac{1}{v^2}\xi\right)\,\ee with
$\delta=vt+\xi^2$. It is easily seen that $\Jac
\left(\tau_i^{(2)}\right) =1$, $i=0,1$. We have checked with MAPLE
that indeed $\tau_i^{(2)}\in \Aut_{K_i} K_i[t,\xi],\,\,i=0,1$,
with the inverse map $\left(\tau_1^{(2)}\right)^{-1}$ given by:
\be\label{inv1} \left(\tau_1^{(2)}\right)^{-1}: (t,\xi)\longmapsto
\left(\frac{1}{v^2}t+x\xi -v^3\xi^2,\,\,\,v^2\xi\right)\,.\ee In
view of (\ref{tf}) and (\ref{inv1}), from
$\left(\tau_0^{(2)}\right)^{-1}= \left(\tau_1^{(2)}\right)^{-1}
\circ\varphi_{10}^{(2)}$
 we obtain:
\be\label{inv0} \left(\tau_0^{(2)}\right)^{-1}
:(t,\xi)\longmapsto \ee
$$\left(x\xi-v^3\xi^2+\frac{2v}{x}t\xi-
\frac{2v^2}{x^3}t^2\xi+
\frac{2}{x^4}t^3-\frac{v}{x^6}t^4,
\quad v^2\xi-\frac{1}{x}t+
\frac{v}{x^3}t^2\right)\,.$$
Thus by Proposition \ref{triv}(a),
$\lambda_n^{(2)}$ is trivial for
any $n\ge 1$, as stated. \qed\medskip

As a corollary of Propositions \ref{equi}
and \ref{trivmain} we recover
the result of V\'en\'ereau \cite{Ve}:

\bcor \label{vcor} For every $n\ge 3$,
the polynomial $v_n$ is an
$x$-variable of the polynomial ring
$\C[x,y,z,u]$. More precisely, the map
$$\alpha_n: \A_\C^4\to \A_\C^4,
\qquad (x,y,z,u)\longmapsto
\left(x,v_n,\zeta^{(n)},\theta^{(n)}\right)\,,$$
is an $x$-automorphism of the ring $\C[x,y,z,u]$,
where
\be\label{theta}
\theta^{(n)}:=u-x^{n-3}t\left(xw^2+
\eta(\zeta^{(n)}-z+x^{n-1}tw)\right)\,.
\ee

\ecor

\proof Indeed, in view of (\ref{pn}),
for $n\ge 3$ we have $\varphi_{10}
=\varphi_{10}^{(2)}$ and so, $\lambda_n=
\lambda_n^{(2)}$ is trivial.

\noindent Furthermore, by virtue of (\ref{inv1}) and Proposition
\ref{loctri}(a) the trivialization $\psi_n$ of $\lambda_n$:
\be\label{etri} \psi_n:=\left(\tau_1^{(2)}\right)^{-1}\circ
\varphi_1^{(n)}=\left(\frac{t+x\zeta^{(n)}
-v_n(\zeta^{(n)})^2}{v_n^2},\,\,\zeta^{(n)}\right)
=:\left(\theta^{(n)},\,\,\zeta^{(n)}\right)\ee extends to a
morphism $$\Psi_n: \A_\C^4 =\Spec \C[x,y,z,u] \to \A_\C^2=\Spec
\C[\theta^{(n)},\,\zeta^{(n)}]$$ providing an automorphism
$\alpha_n:=(\Phi_n,\,\Psi_n)\in \Aut (\A_\C^4)$ (see the proof of
Proposition \ref{triv}). The explicit expression (\ref{theta}) for
$\theta^{(n)}$ was found by applying formulas
(\ref{w})-(\ref{eta}), (\ref{zetan})-(\ref{vn}),  (\ref{zetan''})
and (\ref{etri}); it was checked with MAPLE. \qed
\end{document}